\documentclass[12pt,leqno]{article}
\usepackage{amstex}
\def\newsection{{\setcounter{equation}{0}}\section}

\font\germ=eufm10 scaled \magstep1
\def\goth#1{\hbox{\germ#1}}

\newcommand{\on}{\operatorname}

\def\Gsl{{\goth{sl}}}

\def\Gg{{\goth{g}}}

\def\Gt{{\goth{t}}}

\def\BZ{{\Bbb Z}}
\def\Z{\BZ}
\def\lam{\lambda}
\def\Lam{\Lambda}

\def\lan{\langle}
\def\ran{\rangle}

\def\on{\operatorname}
\def\gge{>\kern-3pt>}

\newcommand{\IM}{{\on{Im}}}

\def\wt{{\on{wt}}}

\def\max{{\on{{max}}}}
\renewcommand{\min}{{\on{min}}}
\newcommand{\cl}{{\on{cl}}}

\def\te{\tilde e}
\def\tf{\tilde f}

\def\beq{\begin{eqnarray}}
\def\beqn{\begin{eqnarray*}}
\def\endeq{\end{eqnarray}}
\def\endeqn{\end{eqnarray*}}
\def\eq{\begin{eqnarray}}
\def\eqn{\begin{eqnarray*}}
\def\nn{\nonumber}
\def\proof{\noindent{\it Proof.}\quad}
\def\qed{\hfill{Q.E.D.}\par\medskip}
\newtheorem{lemma}{Lemma}[section]
\newtheorem{definition}[lemma]{Definition}
\newtheorem{corollary}[lemma]{Corollary}
\newtheorem{sublemma}[lemma]{Sublemma}
\newtheorem{proposition}[lemma]{Proposition}
\newtheorem{theorem}[lemma]{Theorem}
\newtheorem{conjecture}{Conjecture}

\def\Conjecture{\begin{conjecture}}
\def\enconjecture{\end{conjecture}}

\def\Lemma{\begin{lemma}}
\def\enlemma{\end{lemma}}
\def\Definition{\begin{definition}}
\def\Def{\begin{definition}}
\def\endefinition{\end{definition}}
\def\Corollary{\begin{corollary}}
\def\encorollary{\end{corollary}}
\def\Sublemma{\begin{sublemma}}
\def\ensublemma{\end{sublemma}}
\def\Proposition{\begin{proposition}}
\def\Prop{\begin{proposition}}
\def\enproposition{\end{proposition}}
\def\enprop{\end{proposition}}
\def\Theorem{\begin{theorem}}
\def\entheorem{\end{theorem}}

\def\Us{U'_q(\Gg)}

\def\rank{{\rm{rank}}}

\def\semi{\hbox{$\,{\vrule height5.9pt depth.6pt}\kern-1.4pt\times$}}

\def\eps{\varepsilon}

\def\tPsi{\tilde\Psi}

\def\xio{\xi{\raise-2.7pt\hbox{${}_{0}$}}}
\def\refl#1{s{\raise-2.5pt\hbox{${}_{#1}$}}}

\def\hookdownarrow%
{{\raise3pt\hbox{$\cap$}\kern-9pt\raise-4pt\hbox{$\downarrow$}}}
\def\doublevrule
{\raise-5pt\hbox{$\vrule height15pt\hskip2pt\vrule height15pt$}}

\newcommand{\Emax}{E^\max}

\newcommand{\strictsubset}%
{\raise-2pt\hbox{${{\displaystyle\subset}\atop{\not=}}$}}
\newcommand{\varp}{\varphi}
\newcommand{\usl}{U'_q(\widehat\Gsl_2)}
\newcommand{\nid}{\noindent}
\newcommand{\isomo}{\raise2pt\hbox{$\underrightarrow{\,\,\sim\,\,}$}}
\providecommand{\bysame}{\makebox[3em]{\hrulefill}\thinspace}

\begin{document}

\baselineskip=16pt

\title{\bf Quantized Affine Algebras and \\ Crystals with Head}
\author{Seok-Jin Kang\begin{thanks} {Supported in part by 
Basic Science Research Institute 
Program,  Ministry of Education of Korea,  BSRI-97-1414,
and GARC-KOSEF at Seoul National University.}
\end{thanks}
and  Masaki Kashiwara$^{\dag}$
\cr
{$^*$ Department of Mathematics} \cr
{Seoul National University} \cr
{Seoul 151-742, Korea} \cr
{and}\cr
{$^{\dag}$ Research Institute for Mathematical Sciences} \cr
{Kyoto University, Kyoto 606, Japan} }
\date{October 2, 1997}

\maketitle
\renewcommand{\thefootnote}{}

\begin{abstract}
Motivated by the work of Nakayashiki on the inhomogeneous vertex models
of 6-vertex type, we introduce the notion of crystals with head.
We show that the tensor product of the highest weight crystal $B(\lam)$ of
level $k$ and the perfect crystal $B_{l}$ of level $l$ 
is isomorphic to the tensor product of the perfect crystal $B_{l-k}$
of level $l-k$ and the highest weight crystal $B(\lam')$ of level $k$. 

\end{abstract}





\newsection{Introduction}

In \cite{Nak}, Nakayashiki studied the inhomogeneous vertex models of 
$6$-vertex type, and he explained the degeneration of the ground states
from the point of view of the representation theory
as follows.
Let $V(\Lambda_i)$ be the irreducible $\usl$-module
with the highest weight $\Lambda_i$ ($i=0,1$) of level $1$, and
let $V_s$ be the $(s+1)$-dimensional $\usl$-module.
Then there exists an intertwiner
$$\Phi(z):(V_{s-1})_z\otimes V(\Lambda_i)\to V(\Lambda_{i+1})\otimes
(V_s)_z.$$
He identified $(V_{s-1})_z$ with the degeneration of the ground states.

The $q=0$ limit can be described in terms of crystal bases.
Let $B_s$ be the crystal base of $V_s$, and let $B(\Lambda_i)$
be the crystal base of $V(\Lambda_i)$.
Then 
we have an isomorphism of crystals
$$B_{s-1}\otimes B(\Lambda_i)\cong
B(\Lambda_{i+1})\otimes B_s.$$

\vskip 2mm

The purpose of this paper is to generalize the above result
on crystals in a more general situation,
replacing $\usl$ with  quantized affine algebras
$\Us$, $B(\Lambda_i)$ with the crystals of the 
integrable highest weight representations
of arbitrary positive level, 
and $B_s$ with perfect crystals.

\vskip 2mm

The crystal of the 
integrable highest weight representation
has a unique highest weight vector.
Namely, it contains a unique vector $b$ such that
$\te_ib=0$ for all $i$ and
all the other vectors can be obtained from $b$ 
by applying $\tf_i$'s successively.
However, neither  
$B_{s-1}\otimes B(\Lambda_i)$ nor $B(\Lambda_{i+1})\otimes B_s$
has such properties. Instead, they satisfy weaker properties:
the highest weight vector has to be replaced
with a subset consisting of several vectors,
which we call the {\it head}.
This is a combinatorial phenomenon corresponding to
the degeneration of the ground states in the exactly solvable models.

\vskip 2mm

Let $B$ be a crystal.
For $b\in B$, let $E(b)$ be the smallest subset of $B$
containing $b$ and stable under the $\te_i$'s.
We say that $B$ has a {\it head} if $E(b)$ is a finite set for any 
$b\in B$.
For such a crystal, we define its head $H(B)$ to be 
$\{b\in B| \ E(b')=E(b)\ \text{for every $b'\in E(b)$}\}$.
Then the head replaces the role of highest weight vectors:
all the vectors in $B$ can be obtained from vectors 
in the head by applying $\tf_i$'s successively. 

\vskip 2mm

If $D$ is a finite regular crystal
and $B(\lambda)$ is the crystal of the integrable highest weight 
representation with highest weight $\lambda$ of level $k$,
then $D\otimes B(\lambda)$ has a head and
its head is given by $D\otimes u_\lambda$,
where $u_\lambda$ is the highest weight vector of $B(\lambda)$.
However, if we change the order of the tensor product,
the situation is completely different.
The crystal $B(\lambda)\otimes D$ has a head, 
but $u_\lambda\otimes D$ is not the head in general.
In this paper, we prove that, for a perfect crystal $B_{l}$ of
level $l>k$,
$B(\lambda)\otimes B_{l}$
is isomorphic to the crystal $B_{l-k} \otimes B(\lam')$ for 
another dominant integral weight $\lam'$ of level $k$ and 
the perfect crystal $B_{l-k}$ of level $l-k$
(see Theorem \ref{theorem:iso} for more precise statements).
\vskip 2mm

The proof is based on the theory of {\it coherent families of perfect
crystals} developed in \cite{KKM} and the 
characterization of crystals
of the form $D\otimes B(\lambda)$.
We introduce the notion of {\it regular head} (Definition \ref{reg head}),
and we prove that any connected regular crystal with regular head
is isomorphic to a crystal of the form  $H(B)\otimes B(\lambda)$
for some dominant integral weight $\lambda$ 
(see Theorem \ref{Theorem 4.7} for more precise statements).
Then, we check the regularity condition for the coherent families of
perfect crystals.

\vskip 5mm
\noindent
{\it Acknowledgments.} 
The first author would like to express his gratitude to the members
of Research Institute for Mathematical Sciences, Kyoto University for
their hospitality during his stay in the winter and 
the summer of 1997, and the second author thanks the Department of
Mathematics of Seoul National University for their hospitality during
his visit in the fall of 1997.
We would also like to thank T. Miwa for many stimulating discussions.

\vskip 1cm

\newsection{Quantized Affine Algebras}

Let $I$ be a finite index set and $A=(a_{ij})_{i,j\in I}$
a generalized Cartan matrix of affine type. 
We choose a vector space $\Gt$ of dimension $|I|+1$,
and let $\Pi=\{\alpha_i | \ i\in I\}$ and 
$\Pi^{\vee}=\{ h_i | \ i\in I\}$ be linearly independent 
subsets of $\Gt^*$ and $\Gt$, respectively,
satisfying $\langle h_i, \alpha_j \rangle = a_{ij}$
for all $i, j \in I$. The $\alpha_i$ (resp. $h_i$)
are called the {\it simple roots} (resp. {\it simple coroots}), 
and the free abelian group $Q=\bigoplus_{i\in I} \Bbb{Z} \alpha_i$
(resp. $Q^{\vee}=\bigoplus_{i\in I} \Z h_i$) is called the 
{\it root lattice} (resp. {\it dual root lattice}).
We denote by $\delta=\sum_{i\in I} a_i \alpha_i \in Q$  the
smallest positive {\it imaginary root} and $c=\sum_{i\in I}
a_i^{\vee} h_i \in Q^{\vee}$ the {\it canonical central element}
(cf. \cite[Chapter 6] {Kac}).
Set $\Gt^*_{\cl}=\Gt^* / \Bbb{C} \delta$ and let $\cl: \Gt^* \to 
\Gt^*_{\cl}$ be the canonical projection. We denote by 
$\Gt^{*0}=\{ \lambda \in \Gt^* | \ \lan c, \lambda \ran =0 \}$
and $\Gt^{*0}_{\cl}=\cl(\Gt^{*0})$. 

\vskip 2mm

Let $P=\{\lambda \in \Gt^* \ | \ \langle h_i, \lambda \rangle
\in \Z \ \ \text {for all} \ i \in I \}$
be the {\it weight lattice} and 
$P^{\vee}=\{h \in \Gt \ | \ \langle h, \alpha_i \rangle
\in \Z \ \ \text {for all} \ i \in I \}$ 
be the {\it dual weight lattice}.
Note that $\alpha_i, \Lambda_i \in P$ and $h_i \in P^{\vee}$,
where $\Lambda_i \in \Gt^*$ are linear forms satisfying
$\langle h_j, \Lambda_i \rangle =\delta_{ij}$ $(i,j\in I)$. 
Set $P_{\cl}=\cl(P) ={\on{Hom}}(Q^{\vee}, \Bbb {Z})
\subset \Gt^*_{\cl}$, 
$P^0=\{ \lambda \in P | \ \lan c, \lambda \ran =0 \}
\subset \Gt^{*0}$, and $P_{\cl}^0=\cl(P^0) \subset \Gt^{*0}_{\cl}$.

\vskip 2mm

Since the generalized Cartan matrix $A$ is symmetrizable, there
is a non-degenerate symmetric bilinear form $(\ , \ )$ on
$\Gt^*$ satisfying
$$\langle h_i, \lambda \rangle =\frac {2(\alpha_i, \lambda)}
{(\alpha_i, \alpha_i)} \ \ \text {for all} \ i\in I, \ 
\lambda \in \Gt^*.$$
We normalize the bilinear form so that we have 
$$(\delta,\lam)=\lan c,\lam\ran.$$
Note that $\Gt^{*0}_{\cl}$ has a non-degenerate symmetric bilinear form
induced by that on $\Gt^*$. 
We take the smallest positive integer $\gamma$
such that
$\gamma(\alpha_i, \alpha_i) / 2$ is a positive integer for all $i\in I$. 

\vskip 2mm 

\begin{definition} 
\textnormal {
The {\it quantized affine algebra} $U_q({\Gg})$ is the associative 
algebra with 1 over ${\Bbb C}(q^{1/\gamma})$ generated by the elements
$e_i$, $f_i$ $(i\in I)$ and $q(h)$ $(h\in \gamma^{-1} P^{\vee})$
satisfying the following defining relations:
\begin{equation}
\begin{aligned}
& q(0)=1, \ q(h)q(h')=q(h+h') \ \ \ (h,h' \in \gamma^{-1}P^{\vee}), \\
& q(h) e_{i} q(-h) =q^{\langle h, \alpha_i \rangle} e_i, \ \ \\
& q(h) f_{i} q(-h) =q^{-\langle h, \alpha_i \rangle} f_i \ \ \
(h\in\gamma^{-1} P^{\vee}, i\in I), \\
& [e_i, f_j]=\delta_{ij} \frac{t_i-t_i^{-1}}{q_i-q_i^{-1}}
\ \ \ (i,j\in I), \\
& \sum_{k=0}^{1-a_{ij}} (-1)^k e_{i}^{(k)} e_j e_i^{(1-a_{ij}-k)}
=\sum_{k=0}^{1-a_{ij}} (-1)^k f_{i}^{(k)} f_j f_i^{(1-a_{ij}-k)}
=0 \ \ (i\neq j),
\end{aligned}
\end{equation}
where $q_i=q^{(\alpha_i, \alpha_i)/2}$, 
$t_i=q(\frac{(\alpha_i, \alpha_i)}{2} h_{i})$, 
$e_{i}^{(k)}=e_i^k /[k]_i!$, $f_i^{(k)}=f_i^k/[k]_i!$,
$[k]_i =\frac {q_i^{k}-q_i^{-k}}{q_i-q_i^{-1}}$,
and $[k]_i!=[1]_i [2]_i \ldots [k]_i$ for all $i\in I$.}
\end{definition}

The quantized affine algebra $U_q(\Gg)$ has a Hopf algebra structure
with comultiplication $\Delta$, counit $\varepsilon$, and antipode $S$
defined by
\begin{equation}
\begin{aligned}
& \Delta(q(h))=q(h) \otimes q(h),\\
& \Delta (e_i)=e_{i} \otimes t_{i}^{-1} + 1 \otimes e_{i}, \ \ 
\Delta (f_i)=f_{i} \otimes 1 + t_{i} \otimes f_{i},\\
& \varepsilon(q(h))=1, \ \ \varepsilon(e_i)=\varepsilon(f_i)=0, \\
&  S(q(h))=q(-h), \ \ S(e_i)=-e_i t_i, \ \ S(f_i)=-t_i^{-1} f_i
\end{aligned}
\end{equation}
for all $h\in\gamma^{-1} P^{\vee}$, $i\in I$. 

\vskip 2mm

We denote by $U_q'(\Gg)$ the subalgebra of $U_q(\Gg)$
generated by $e_i$, $f_i$ $(i\in I)$ and $q(h)$ $(h\in\gamma^{-1} Q^{\vee})$,
which will also be called the {\it quantized affine algebra}.

\vskip 2mm

A $U'_q(\Gg)$-module $M$ is called {\it integrable} if it has 
the {\it weight space decomposition} 
$M=\bigoplus_{\lambda \in P_{\cl}} M_{\lambda}$, where
$M_{\lambda}=\{ u\in M| \ q(h) u = q^{\lan h, \lam \ran} u
\ \ \text {for all} \ h\in\gamma^{-1}Q^{\vee} \}$,
and $M$ is $U'_q(\Gg)_{i}$-locally finite (i.e., 
$\dim U'_q(\Gg)_{i} u < \infty$ for all $u\in M$) for
all $i\in I$, where $U'_q(\Gg)_{i}$ denotes the subalgebra
of $U'_q(\Gg)$ generated by $e_i$, $f_i$, and $t_i$.

\vskip 2cm

\newsection{Crystals with Head}

In studying the structure of integrable representations of quantized 
affine algebras, the {\it crystal base theory} developed in \cite {Kas1}
provides a very powerful combinatorial method. In this section,
we develop the theory of {\it crystals with head}. 
We first recall the definition of {\it crystals} given in \cite{Kas2}.

\begin{definition} 
\textnormal {A {\it crystal} $B$ is a set together with the maps
$\wt: B \rightarrow P$, 
$\varepsilon_{i}: B \rightarrow {\Bbb Z} \sqcup \{-\infty \}$,
$\varphi_{i}: B \rightarrow {\Bbb Z} \sqcup \{-\infty \}$,
${\tilde e_{i}}: B \rightarrow B \sqcup \{0\}$, 
${\tilde f_{i}}: B \rightarrow B \sqcup \{0\}$ $(i\in I)$
satisfying the axioms:
\begin{equation}
\begin{aligned}
& \lan h_i, \wt(b) \ran =\varphi_{i}(b) -\varepsilon_{i}(b) 
\ \ \text {for all} \ b\in B, \\
& \wt({\tilde e_{i}}b)=\wt(b)+\alpha_{i} \ \ \text {for} \ 
b\in B \ \ \text {with} \ {\tilde e_{i}} b \in B, \\
& \wt({\tilde f_{i}}b)=\wt(b)-\alpha_{i} \ \ \text {for} \ 
b\in B \ \ \text {with} \ {\tilde f_{i}} b \in B, \\
& {\tilde f_{i}}b=b' \ \ \text {if and only if} \ \ 
b=\tilde e_{i} b' \ \ \text {for} \ b, b'\in B, \\
& \tilde e_i b=\tilde f_i b=0 \ \ \text {if} \ 
\ \varepsilon_{i}(b)=-\infty.
\end{aligned}
\end{equation}}
\end{definition}

\vskip 2mm

\begin{definition} 
\textnormal {For two crystals $B_1$ and $B_2$, a {\it morphism} 
of crystals from $B_1$ to $B_2$ is a map 
$\psi: B_1 \sqcup \{0\} \to B_2 \sqcup \{0\}$ such that
\begin{equation}
\begin{aligned}
& \psi(0)=0, \\
& \psi ({\tilde e_{i}}b)=\tilde e_i \psi(b) \ \ \text {for} \ 
b\in B_1 \ \ \text {with} \ \ {\tilde e_{i}} b \in B_1, \ 
\psi(b) \in B_2, \ \psi(\tilde e_i b) \in B_2, \\
& \psi ({\tilde f_{i}}b)=\tilde f_i \psi(b) \ \ \text {for} \ 
b\in B_1 \ \ \text {with} \ \ {\tilde f_{i}} b \in B_1, \ 
\psi(b) \in B_2, \ \psi(\tilde f_i b) \in B_2, \\
& \wt(\psi(b))=\wt(b) \ \ \text {for} \ 
b\in B_1 \ \ \text {with} \ \ \psi(b) \in B_2, \\
& \varepsilon_i(\psi(b))=\varepsilon_i(b), 
\ \ \varphi_i(\psi(b))=\varphi_i(b) \ \ \text {for} \ 
b\in B_1 \ \ \text {with} \ \ \psi(b) \in B_2. \\
\end{aligned}
\end{equation}}
\end{definition}
A morphism $\psi: B_1 \to B_2$ is called an {\it embedding} if 
the map $\psi:B_1 \sqcup \{0\} \to B_2 \sqcup \{0\}$ is injective. 
In this case, we call $B_1$ a {\it subcrystal} 
of $B_2$.

\vskip 2mm

For two crystals $B_1$ and $B_2$, we define their {\it tensor product}
$B_1 \otimes B_2$ as follows. The underlying set is $B_1 \times B_2$.
For $b_1 \in B_1$, $b_2 \in B_2$, we write $b_1 \otimes b_2$ for
$(b_1, b_2)$ and we understand $b_1 \otimes 0 = 0 \otimes b_2 =0$.
We define the maps
$\wt: B_1 \otimes B_2  \rightarrow P$, 
$\varepsilon_{i}: B_1 \otimes B_2 \rightarrow {\Bbb Z} \sqcup \{-\infty \}$,
$\varphi_{i}: B_1 \otimes B_2 \rightarrow {\Bbb Z} \sqcup \{-\infty \}$,
${\tilde e_{i}}: B_1 \otimes B_2 \rightarrow 
B_1 \otimes B_2 \sqcup \{0\}$, 
${\tilde f_{i}}: B_1 \otimes B_2 \rightarrow 
B_1 \otimes B_2 \sqcup \{0\}$ $(i\in I)$ as follows:
\begin{equation}
\begin{aligned}
& \wt(b_1 \otimes b_2)  =\wt(b_1) + \wt(b_2), \\
& \varepsilon_i (b_1 \otimes b_2)= \max (\varepsilon_i (b_1),
\ \varepsilon_i(b_2)-\lan h_i, \wt(b_1) \ran ), \\
& \varphi_i (b_1 \otimes b_2)= \max (\varphi_i (b_2),
\ \varphi_i(b_1) + \lan h_i, \wt(b_2) \ran ), \\
& \tilde e_{i}(b_1 \otimes b_2)=\cases \tilde e_i b_1 \otimes b_2
\ \ & \text {if} \ \varphi_i(b_1) \ge \varepsilon_i (b_2), \\
b_1 \otimes \tilde e_i b_2 \ \ & \text {if} \ 
\varphi_i(b_1) < \varepsilon_i (b_2), \endcases \\
& \tilde f_{i}(b_1 \otimes b_2)=\cases \tilde f_i b_1 \otimes b_2
\ \ & \text {if} \ \varphi_i(b_1) > \varepsilon_i (b_2), \\
b_1 \otimes \tilde f_i b_2 \ \ & \text {if} \ 
\varphi_i(b_1) \le \varepsilon_i (b_2). \endcases 
\end{aligned}
\end{equation}

\vskip 2mm 
In the sequel, we will 
only consider the crystals over the quantized affine algebra 
$U_q'(\Gg)$. Hence the weights of crystals will be elements of $P_{\cl}$.
For example, for $\lam \in P_{\cl}$, consider the
set $T_{\lam}=\{t_{\lam} \}$ with one element. 
Define $\wt(t_{\lam})=\lam$, $\varepsilon_i (t_{\lam})
=\varphi_i(t_{\lam})=-\infty$, and $\tilde e_i (t_{\lam})
=\tilde f_i (t_{\lam})=0$ $(i\in I)$. Then $T_{\lam}$ 
is a crystal and we have $T_{\lam} \otimes T_{\lam'}
\cong T_{\lam+\lam'}$. 

\vskip 2mm

For a dominant integral weight $\lam$,
we denote by $B(\lam)$ the crystal associated with
the integrable highest weight 
representation with highest weight $\lambda$,
and $u_\lam$ the highest weight vector of $B(\lam)$.
The highest weight vector $u_{\lam}$ is the unique element of $B(\lam)$
with weight $\lam$ satisfying $\tilde e_{i} u_{\lam}=0$ for all $i\in I$. 
\vskip 2mm

For a subset $J$ of $I$, 
we denote by $U'_q(\Gg_J)$ the subalgebra of $\Us$ generated by
$e_i$, $f_i$, and $t_i$ ($i\in J$).
Note that if $J\strictsubset I$, then $\Gg_J$ 
is a finite-dimensional semisimple Lie algebra.
Similarly, for a subset $J$ of $I$, we denote by $B_{J}$
the crystal $B$ equipped with the maps $\wt$, $\varepsilon_{i}$,
$\varphi_{i}$, $\tilde e_{i}$, and $\tilde f_{i}$ for
$i\in J$. 
We say that a crystal $B$ over $\Us$
is {\it regular}
if, for any $J\strictsubset I$,
$B_{J}$ is isomorphic to
the crystal associated with
an integrable $U'_q(\Gg_J)$-module.
This condition is equivalent to saying that
the same assertion holds for any $J\strictsubset I$
with one or two elements (see \cite[Proposition 2.4.4]{KMN1}).

\vskip 2mm 
Let $B$ be a regular crystal.
For $b\in B$, let ${\tilde e_i}^{\max} b={\tilde e_{i}}^k b$ such 
that ${\tilde e_{i}}^k b \neq 0$, ${\tilde e_{i}}^{k+1} b=0$,
and define 
\begin{equation}
\begin{aligned}
E(b)&= \{\te_{i_1}\ldots\te_{i_l}b \ | \ 
\hbox{$l\ge0$ and $i_1,\ldots,i_l\in I$}\}\setminus\{0\},\\
\Emax(b)&=\{\te^\max_{i_1}\ldots\te^\max_{i_l}b \ | \ 
\hbox{$l\ge0$ and $i_1,\ldots,i_l\in I$}\}.
\end{aligned}
\end{equation}
It follows that 
\begin{equation}
\begin{aligned}
 \Emax(b) & \subset E(b), \\
 E(b') &  \subset E(b) \ \ \text {for all} \ b' \in E(b), \\
 \Emax(b') &  \subset \Emax(b) \ \ \text {for all} \ b' \in \Emax(b).
\end{aligned}
\end{equation}

Recall that the Weyl group $W$ acts on the regular crystals (\cite{Kas2}).
For each $i\in I$, the simple reflection $s_i$ acts on the regular 
crystal $B$  by
\begin{equation} \label {3.5}
S_{i}(b)= \cases \tilde f_{i}^{\langle h_i, \wt(b)\rangle} b
\ \ & \text {if} \ \langle h_i, \wt(b) \rangle \ge 0, \\
\tilde e_{i}^{-\langle h_i, \wt(b) \rangle} b \ \ & \text {if}
\ \langle h_i, \wt(b) \ran \le 0. \endcases
\end{equation}
For $w=s_{i_r} s_{i_{r-1}} \ldots s_{i_1} \in W$, its action is given by
$S_{w}=S_{i_r}S_{i_{r-1}} \ldots S_{i_1}$.

\vskip 2mm 

We first prove:

\Lemma \label {Lemma 3.2}
Let $B$ be a finite regular crystal.

{\rm (a)} We have $E(S_{w}(b))=E(b)$ for all $b\in B$, $w\in W$.

\vskip 2mm

{\rm (b)} $E(b)$ is a connected component of $B$
for any $b\in B$.
\enlemma

\proof
(a) It suffices to show that $S_{i}(b) \in E(b)$ for all $b\in B$,
$i\in I$. If $\lambda=\wt(b)$ satisfies $\lan h_i, \lam \ran \le 0$,
then by (\ref{3.5}), our assertion is obvious.
If $ \lan h_i, \lam \ran >0$, take $w=s_{i_l}\ldots s_{i_1}\in W$ 
such that $\lan h_{i_k},s_{i_{k-1}}\ldots s_{i_1} \lam \ran <0$
for $k=1,\ldots,l$ and $s_{i} \lam =w \lam$ (see \cite[Lemma 1.4]{Akas}).
Then for each $n \ge 1$, we have $S_{w} (S_{i}S_{w})^n b \in E(b)$.
Since $S_{i}S_{w}$ has finite order, there exists $n>0$ such that
$(S_{i}S_{w})^n b=b$. Hence 
$S_{i}b = S_{w} (S_{i}S_{w})^{n-1} b \in E(b)$. 

\vskip 2mm

(b) Note that for any $b\in B$, we have
$\tilde f_{i} b = \tilde e_{i}^{\varphi_{i}(b)-1} S_{i} (\te^{\max}_{i} b)$.
By (a), this implies $E(b)$ is stable under $\tilde f_{i}$ for all
$i\in I$. Hence we have the desired result.
\qed

\vskip 2mm 

\begin{definition} 
\textnormal {
We say that a regular crystal 
$B$ has a {\it head} if $E(b)$ is a finite set
for any $b\in B$.
In this case, we define the {\it head} $H(B)$ of $B$ to be 
\eq
&&H(B)=\{b\in B| \ E(b')=E(b)\quad\hbox{for every $b'\in E(b)$}\},
\endeq
and $B$ is called a {\it crystal with head}.}
\end{definition}

\vskip 2mm 
In the following, we prove some of the basic properties of the crystals
with head.

\Lemma 

Suppose that $B$ has a head $H(B)$. 

{\rm (a)} The head $H(B)$ is stable under $\tilde e_{i}$'s $(i\in I)$.

\vskip 2mm

{\rm (b)} $E(b)\cap H(B)\not=\emptyset$ for all $b\in B$.

\vskip 2mm

{\rm (c)} If $b \in H(B)$, then either $\tilde e_{i} b=0$ for
all $i\in I$ or there exist $i_{1}, \ldots, i_{l} \in I$ 
$(l\ge 1)$ such that $b=\tilde e_{i_l} \ldots \tilde e_{i_1}b$.
\enlemma

\proof
(a) If $b\in H(B)$, then $E(b) \subset H(B)$, since for $b'\in E(b)$
and $b''\in E(b') \subset E(b)$, we have $E(b'')=E(b)=E(b')$.

\vskip 2mm

(b) For any $b\in B$, take $b'\in E(b)$ 
such that $E(b')$ has the smallest cardinality.
Then, since $E(b'') \subset E(b')$ for any $b''\in E(b') \subset E(b)$, 
we have $E(b'')=E(b')$, which implies $b'$ belongs to $H(B)$.

\vskip 2mm

(c) If $b\in H(B)$ and $\tilde e_{i_1} b \neq 0$ for some
$i_1 \in I$, then by definition we have $E(b)=E(\tilde e_{i_1}b)$.
Then $b\in E(b)$ implies $b=\tilde e_{i_l} \ldots \tilde e_{i_1}b$
for some $i_2, \ldots, i_l \in I$. 
\qed

\Lemma \label {Lemma 3.5}
Let $B$ be a regular crystal with head and $H$  a  
subset of $B$.

{\rm (a)} If $H$ is stable under $\te_i$'s $(i\in I)$
and $E(b) \cap H \neq \emptyset$ for any $b\in B$, then
$H(B)$ is contained in $H$. 

\vskip 2mm

{\rm (b)} If, in addition, $E(b)=E(b')$ for any $b\in H$
and $b'\in E(b)$, then $H=H(B)$. 
\enlemma

\proof (a) If $b\in H(B)$, take $b'\in E(b)\cap H$.
Then $b\in E(b)=E(b') \subset H$. 

\vskip 2mm

(b) If $b\in H$ and $b'\in E(b) \cap H(B)$, then 
$b\in E(b)=E(b') \subset H(B)$. 
\qed

\Corollary
If $B$ is a finite regular crystal, then $H(B)=B$.
\encorollary

\proof 
We may assume that $B$ is connected. By Lemma \ref{Lemma 3.2}, 
we have $E(b)=B$ 
for all $b\in B$. Hence $H(B)=B$. \qed

\newsection{Structure of Crystals with Head}

Let $\psi: H(B) \hookrightarrow B$ denote the inclusion map. 
\Definition\label{reg head}
\textnormal{
We say that $B$ has a {\it regular head} if 
the head $H(B)$ of $B$ becomes a regular crystal
with the maps $\wt$, $\varepsilon_{i}$,
$\varphi_{i}$, $\tilde e_{i}$, $\tilde f_{i}$ $(i\in I)$
defined by 
\begin{equation}
\begin{aligned}
\tilde e_i b& =\psi^{-1}(\tilde e_i \psi(b)), \\
 \tilde f_i b&=\cases  \psi^{-1}(\tilde f_i \psi(b)) \ \ 
 & \text {if} \ \tilde f_{i} \psi(b) \in H(B), \\
 0  \ \ & \text {otherwise}, \endcases \\
 \varepsilon_{i}(b)&=\varepsilon_{i}(\psi(b)), \\
 \varphi_{i}(b)&=\max \{ k\ge 0| \ \tilde f_{i}^{k} b \in H(B) \} \\
&=\max \{ k \ge 0 | \ 
b\in \tilde e_{i}^k H(B) \}, \\ 
 \wt(b)&=\sum_i(\varphi_i(b) -\varepsilon_{i}(b))\Lam_i\in P_\cl.
\end{aligned}
\end{equation}
}
\end{definition}

Let $b\in H(B)$. Then $E(b) \subset H(B)$.
If $b'\in E(b)$ satisfies $\tf_i b'\in H(B)$, then 
$\tf_i b' \in E(\tf_i b')=E(\te_i \tf_i b')=E(b')=E(b)$,
and hence $E(b)$ is stable under $\tf_i$'s $(i\in I)$.
Therefore the connected components of $H(B)$ are of the form $E(b)$. 

\vskip 2mm

Let $E(b_{0})$ be a connected component of $H(B)$
and set $W(b)=\wt(\psi(b))-\wt(b)$ for $b\in E(b_{0})$,
where $\psi:H(B) \hookrightarrow B$ is the inclusion map.
Note that, for all $i,j\in I$, we have
\begin{equation}
\begin{aligned}
\lan h_i, W(\te_{j} b) \ran &
=\lan h_i, \wt(\psi(\te_j b)) \ran -\lan h_i,  \wt(\te_j b) \ran \\
&=\lan h_i, \wt(\psi(b)) + \alpha_j \ran
-\lan h_i,  \wt(b) +\alpha_j \ran \\
&=\lan h_i, \wt(\psi(b))-\wt(b) \ran 
=\lan h_i, W(b) \ran.
\end{aligned}
\notag
\end{equation}
Hence, $W(\te_j b)=W(b)$ for all $j\in I$, which implies
$W(b)$ is constant on $E(b_{0})$. 

\vskip 2mm
Let $\lambda_0=\wt(\psi(b_{0}))-\wt(b_{0})$.
Since 
\begin{equation}
\begin{aligned}
\lan h_i, \lambda_0 \ran &
=\lan h_i, \wt(\psi(b_{0})) - \wt(b_{0}) \ran \\
&=\varphi_{i} (\psi(b_{0}))-\varphi_{i}(b_{0}) \ge 0,
\end{aligned}
\notag
\end{equation}
$\lambda_0$ is dominant integral. 
We will show that there exists a unique 
embedding of regular crystals $E(b_0) \otimes B(\lambda_0)
\rightarrow B$ sending $b\otimes u_{\lambda_{0}}$
to $\psi(b)$ for all $b\in E(b_0)$, where $u_{\lambda_0}$ is the
highest weight vector of $B(\lambda_0)$.

\vskip 2mm 

Let $D$ be a finite regular crystal,
and let $\lam$ be a dominant integral weight.
We denote by $B(\lam)$ the crystal associated with
the integrable highest weight $\Us$-module $V(\lam)$ 
with highest weight $\lambda$, 
and let $u_{\lam}$ be the highest weight vector of
$B(\lam)$. 

\Lemma \label{Lemma 4.2}
For any $b\in D\otimes B(\lam)$, we have
$$\Emax(b)\cap (D\otimes u_\lam) \not=\emptyset.$$
\enlemma
\proof
If it were not true, there would exist 
$b=b_1\otimes b_2\in D\otimes B(\lam)$ such that
$\Emax(b)\subset D\otimes b_2$ and $b_2 \in B(\lam)\setminus\{u_\lam\}$.
By the tensor product rule, this implies
$\Emax(b)=\Emax(b_1)\otimes b_2$.
Since $b_2 \not=u_\lam$, there exists $i\in I$ such that
$\eps_i(b_2)>0$. 
Take $b' \in\Emax(b_1)$ such that $\varp_i(b')=0$. 
Such a $b'$ exists by \cite[Lemma 1.5]{Akas}.
Then we have
$\te_i(b'\otimes b_2)=b' \otimes (\te_i b_2)$,
which contradicts $\te^\max_i(b' \otimes b_2)\in D\otimes b_2$.
\qed

\Lemma
The regular crystal $D\otimes B(\lam)$ has a regular head and
$H(D\otimes B(\lam))=D\otimes u_\lam$, which is isomorphic to $D$
as a crystal. 
\enlemma

\proof
Since $E(b_1\otimes b_2)\subset E(b_1)\otimes E(b_2)$,
$D\otimes B(\lam)$ has a head.
The second assertion follows from Lemma \ref{Lemma 3.2}, 
Lemma \ref{Lemma 3.5} and
Lemma \ref{Lemma 4.2}. Hence $D\otimes B(\lam)$ has a regular head.
\qed

\Prop \label {Prop 4.4}
Let $D$ be a finite regular crystal,
and let $\lam$ be a dominant integral weight.
Then for every
$b\in D\otimes B(\lam)$,
there exists a positive integer $N$ such that
$\te^\max_{i_N}\ldots\te^\max_{i_1}b \in D\otimes u_\lam$
if  $\te^\max_{i_k}\ldots\te^\max_{i_1}b\not=
\te^\max_{i_{k-1}}\ldots\te^\max_{i_1}b$ for $1\le k \le N$.
\enprop

\proof
If the proposition were false,
there would exist 
$b\in (D\otimes B(\lam))\setminus (D\otimes u_\lam)$
and $l>0$ such that
\eq
&&\text{$b=\te^\max_{i_l}\ldots\te^\max_{i_1}b \quad$
 and $\quad\te^\max_{i_k}\ldots\te^\max_{i_1}b\not=
\te^\max_{i_{k-1}}\ldots\te^\max_{i_1}b$ \ for \ $k =1,\ldots,l$.}
\endeq
Set $\te^\max_{i_k}\ldots\te^\max_{i_1}b=b_k\otimes b'$
with $b_k\in D$ and $b'\in B(\lam)$.
Then $b'$ does not depend on $k$
and we have
$b_k=\te^\max_{i_k}b_{k-1}$.
Since $D$ is a finite crystal, all of its weights have level $0$.
Hence the square lengths of its weights are well-defined.

\vskip 2mm
Since $\wt(b_{k})=\wt(b_{k-1})+\varepsilon_{i_k} (b_{k-1}) \alpha_{i_k}$,
we have
\eq \label {4.3}
(\wt(b_{k}), \wt(b_{k}))&=& (\wt(b_{k-1}), \wt(b_{k-1})) 
+ 2 \varepsilon_{i_k} (b_{k-1}) (\wt (b_{k-1}), \alpha_{i_k}) \\
&&+ \varepsilon_{i_k} (b_{k-1})^2 (\alpha_{i_k}, \alpha_{i_k}) \nn\\
&=& (\wt(b_{k-1}), \wt(b_{k-1})) 
+ \varepsilon_{i_k}(b_{k-1}) (\alpha_{i_k}, \alpha_{i_k}) 
\lan h_{i_k}, \wt(b_{k-1}) \ran \nn \\
&& +\varepsilon_{i_k} (b_{k-1})^2 (\alpha_{i_k}, \alpha_{i_k})\nn \\
&= &(\wt(b_{k-1}), \wt(b_{k-1})) 
+ (\alpha_{i_k}, \alpha_{i_k}) \varepsilon_{i_k}(b_{k-1})
 \varphi_{i_k} (b_{k-1})\nn \\
& \ge& (\wt(b_{k-1}), \wt(b_{k-1})) \nn
\endeq
for all  $k\ge 1$.
Hence $(\wt(b_k),\wt(b_k))$ are the same for all $k\ge 1$.
Since (\ref{4.3}) is the equality and $\eps_{i_k}(b_{k-1})>0$,
we have $\varphi_{i_k} (b_{k-1})=0$.
Since $\tilde e_{i_k} (b_{k-1} \otimes b')
=\tilde e_{i_k} b_{k-1} \otimes b'$,
we have $\varphi_{i_k} (b_{k-1}) \ge \varepsilon_{i_k} (b')$,
and hence $\varepsilon_{i_k} (b')=0$.
Write $\wt (\te^\max_{i_l}\ldots\te^\max_{i_1}b)
=\cl(t_{1} \alpha_{i_1} + \cdots + t_{l} \alpha_{i_l})+\wt(b)$.
Since $\wt(b)=\wt(\te^\max_{i_l}\ldots\te^\max_{i_1}b)$, 
$t_{1} \alpha_{i_1} + \cdots + t_{l} \alpha_{i_l}$
is a multiple of the null root $\delta$, 
which implies $\{i_1,\ldots,i_l\}=I$.
Hence $\varepsilon_{i} (b')=0$ for all $i\in I$, which 
contradicts $b' \neq u_{\lam}$.
\qed

\vskip 2mm
Note that the subcrystal $D\otimes u_\lam$ of $D \otimes B(\lam)$ 
is isomorphic to the crystal 
$D\otimes T_\lam$, where $T_{\lam}$ denotes the crystal with
a single element $t_{\lam}$ of weight $\lam$ and with
$\varepsilon_i(t_{\lam})=\varphi_i(t_{\lam})=-\infty$. 
Let $B$ be a regular crystal. 
In the next theorem, we will show that any morphism
of crystals $\Psi: D \otimes u_{\lam} \rightarrow B$
commuting with the $\tilde e_{i}$'s $(i\in I)$ can be
extended uniquely to a morphism of regular crystals from
$D \otimes B(\lam) \to B$. 

\vskip 2mm

\Theorem \label {Theorem 4.5}
Let $D$ be a finite regular crystal, $B$  a regular crystal, and
$\lam$ a dominant integral weight. Suppose that there is a 
morphism of crystals 
\eqn
&&\Psi:D\otimes u_\lam\to B
\endeqn
such that $\Psi(D \otimes u_{\lam}) \subset B$
and $\Psi$ commutes with the $\te_i$'s $(i\in I)$.

\vskip 2mm
Then, if $\text {rank} \ \Gg >2$, the map $\Psi$ can be 
uniquely extended to a morphism of regular crystals
$$ \tilde \Psi : D\otimes B(\lam) \to B. $$
\entheorem
\newcounter{eqi}
\newcounter{eqii}
\proof
Let $\Sigma$ be the set of pairs $(S,\tilde\Psi)$ satisfying 
the following properties:
\eq
&&D\otimes u_\lam\subset S\subset D\otimes B(\lam),\label{p1}\\
&&\hbox{$\te^\max_iS\subset S$ for any $i \in I$,}\label{p2}\\
&&\hbox{$\tilde\Psi$ is a map from $S$ to $B$
such that $\tilde\Psi|_{D\otimes u_\lam}=\Psi$,}\label{p4}
\setcounter{eqi}{\value{equation}}
\\
&&\hbox{$\wt(\tPsi(b))=\wt(b)$ and $\eps_i(\tPsi(b))=\eps_i(b)$
for any $b\in S$ and $i \in I$,}\label{p5}\\
&&\tPsi(\te^\max_ib)=\te^\max_i\tPsi(b)
\quad\hbox{for any $b\in S$ and $i\in I$.}\label{p6}
\setcounter{eqii}{\value{equation}}
\endeq
Since $\Sigma$ is inductively ordered, 
by Zorn's Lemma, it has a maximal element.
Let $(S,\tilde\Psi)$ be a maximal element.
It is enough to prove that $S$ is the same as $D\otimes B(\lam)$.
Assume that they are different.

\vskip 2mm

First we shall prove 
that there exists
$b\in D\otimes B(\lam)\setminus S$
such that $\te^\max_i(b)\in S\cup\{b\}$ for any $i\in I$.
If it were not true, for any $b\in D\otimes B(\lambda) \setminus S$,
there would exist $i$ such that $\te^\max_i(b)\not\in S\cup\{b\}$.
Let us take $b_{0} \in D \otimes B(\lam) \setminus S$.
Then there is $i_0$ such that $b_1=\te^\max_{i_0}(b)\not\in S\cup\{b_0\}$.
Repeating this we can find a sequence $\{b_k\}$ and $\{i_k\}$
such that $b_{k+1}=\te^\max_{i_k}(b_k)\not\in S\cup\{b_k\}$.
This contradicts Proposition \ref{Prop 4.4}.
Hence there exists 
$b \notin S$ and $\te^\max_{i} b \in S \cup \{b\}$ for all 
$i\in I$. 
We shall choose such a $b$.
\vskip 2mm 

Next we shall show that there exists $i_0$ such that
$\te^\max_{i_0}(b)\in S$.
Assuming the contrary, we shall deduce a contradiction.
Write $b=b_{1} \otimes b_{2}$. If $\te^\max_i b =b$ for
all $i\in I$, then 
$0=\eps_i(b_1\otimes b_2)=\max(\eps_i(b_1),\eps_i(b_2)-\lan h_i,\wt(b_1)\ran)$.
This implies $\eps_i(b_1)=0$ for every $i$, and therefore
$\lan h_i,\wt(b_1)\ran=\varphi_i(b_1)\ge0$.
Since $\lan c,\wt(b_1)\ran=0$, we have
$\lan h_i,\wt(b_1)\ran=0$ for every $i$.
Thus we obtain $\eps_i(b_2)=0$ for every $i$ and hence $b_2=u_\lambda$.
This contradicts $D\otimes u_\lambda\subset S$.

\vskip 2mm

Note that $\varp_{i_0}(\tPsi(\te^\max_{i_0}b))=\varp_{i_0}(\te^\max_{i_0}b)
\ge\eps_{i_0}(b)$.
We define $\tPsi(b)$ to be  
$\tf_{i_0}^{\eps_{i_0}(b)}\tPsi(\te^\max_{i_0}b)\in B$.
We will show that $(S\cup\{b\},\tPsi)$ satisfies
(\ref{p1}--\theeqii).
The properties (\ref{p1}--\theeqi) are automatically satisfied.
For (\ref{p5}), note that
\eqn
&&\wt(\tPsi(b))=\wt(\tPsi(\te^\max_{i_0}b))-\eps_{i_0}(b)\alpha_{i_0}
=\wt(\te^\max_{i_0}b)-\eps_{i_0}(b)\alpha_{i_0}=\wt(b).
\endeqn
We shall show $\eps_i(\tPsi(b))=\eps_i(b)$ for $i\in I$.
Set $J=\{i,i_0\}\strictsubset I$.
Let $K$ be the connected component of $D \otimes B(\lam)$ 
as a $U'_q(\Gg_{J})$-crystal containing $b$.
Then $K$ is a finite set. 
Take a highest weight vector $b_1\in K \subset D\otimes B(\lam)$.
Then since $\te^\max_{i_0} (b) \in S$ and 
$\te^\max_{i} S \subset S$ for all $i\in I$,
$b_1$ lies in $S$.
By (\ref{p5}), $\tPsi(b_1)$ is also a highest weight vector
with respect to the $J$-colored arrows,  
and $\wt(\tPsi(b_1))=\wt(b_1)$. 
Hence the map $b_1\mapsto \tPsi(b_1)$ extends to a morphism of 
$U'_q(\Gg_J)$-crystals $\psi:K\to B$.
Evidently, $\psi|{}_{K\cap S}=\tPsi|{}_{K\cap S}$.
Since 
$$\tPsi(b)=\tilde f_{i_0}^{\varepsilon_{i_0} (b)}
\tPsi(\te^\max_{i_0}b)
=\tilde f_{i_0}^{\varepsilon_{i_0} (b)}
\psi(\te^\max_{i_0}b)=\psi( \tilde f_{i_0}^{\varepsilon_{i_0}(b)}
\te^\max_{i_0}b)=\psi(b),$$
we have the desired property $\eps_i(\tPsi(b))=\eps_i(b)$.

\vskip 2mm 

Finally, let us prove (\ref{p6}). 
If $\te^\max_i(b)\in S$, then
\eqn
&&
\te^\max_i\tPsi(b)=
\te^\max_i\psi(b)=\psi(\te^\max_i(b))=\tPsi(\te^\max_ib).
\endeqn
If $\te^\max_i(b)=b$, then $\varepsilon_i(b)=0$, and hence
$\varepsilon_i(\tPsi(b))=0$. 
Thus $\tPsi(\te^\max_i b)=\tPsi(b)=\te^\max_i \tPsi(b)$.

\qed

\Corollary
Let $B$ be a regular crystal with regular head.
For an arbitrary connected component $E(b_{0})$ of $H(B)$, 
let $\psi: E(b_0) \hookrightarrow H(B)$ be the inclusion map. 
Then there exists a unique embedding of regular crystals 
$\Psi: E(b_{0}) \otimes B(\lam_0) \rightarrow B$
such that $\Psi(b\otimes u_{\lam_0})=\psi(b)$
for any $b\in E(b_0)$.
\encorollary
\proof
Since $E(b_0)$ is finite, the existence and the uniqueness of $\Psi$
follow immediately from Theorem \ref{Theorem 4.5}.
We can also see that $\Psi$ is an embedding by Lemma \ref{Lemma 4.2}.

\qed

\vskip 2mm 

The following theorem describes completely the structure of the regular
crystals with regular head. 

\vskip 2mm

\Theorem \label{Theorem 4.7}  Suppose $\rank \ \Gg>2$. Then any 
regular crystal $B$ with regular head
has the following decomposition:
$$B \cong \bigsqcup_{D} D \otimes B(\lam_D),$$
where $D$ ranges over the connected components of $H(B)$
and $\lambda_D$ is a dominant integral weight.
\entheorem
\proof It suffices to prove that $E(b_0) \otimes B(\lam)$
is connected for all $b_0 \in H(B)$. 
This follows from the fact that $H(E(b_{0}) \otimes B(\lam))
\cong E(b_{0}) \otimes u_{\lam}$ and 
$E(b \otimes u_{\lam})=E(b) \otimes u_{\lam} \ni b_{0} \otimes u_{\lam}$
for any $b\in E(b_{0})$.
\qed

\newsection{Highest Weight Crystals and Perfect Crystals}

Let $k$, $l$ be positive integers, $\lam$ a dominant integral weight
of level $k$, and $B_l$  a {\it perfect crystal} of level $l$.
The definition and the relevant theory of perfect crystals 
can be found in \cite {KKM}, \cite{KMN1} and \cite{KMN2}. 
Consider the
tensor product of regular crystals 
$B(\lambda) \otimes B_l$, where $B(\lam)$ is the crystal for the 
integrable highest weight module $V(\lam)$ over $\Us$ with a 
dominant integral highest weight $\lam$.
If $k\ge l$, it is known that $B(\lam) \otimes B_l$ 
decomposes into a disjoint union of crystals $B(\mu)$, where $\mu$ 
is a dominant integral weight of level $k$. 
In fact, $H(B(\lam) \otimes B_{l})$ is a discrete 
crystal in this case, and coincides with $u_{\lam} \otimes
B_{l}^{\le \lam}$, where 
$B_{l}^{\le \lam}=\{ b\in B_{l} | \ \varepsilon_{i} (b) 
\le \lan h_i, \lam \ran \ \ \text {for all} \ i\in I \}$. 
Hence we have
$$B(\lam) \otimes B_{l} \cong \bigoplus_{b\in B_{l}^{\le \lam}}
B(\lam + \wt (b)).$$
See \cite{KMN1} and \cite{KMN2} for details.

\vskip 2mm

In this work, we will concentrate on the case when $k<l$. 
We first observe:

\Proposition
The crystal $B(\lam) \otimes B_{l}$ has a head and 
$$H(B(\lam) \otimes B_{l}) \subset u_{\lam} \otimes B_{l}.$$
\enproposition

\proof For any $b_{1} \otimes b_{2} \in B(\lam) \otimes B_{l}$,
we have $E(b_{1} \otimes b_{2}) \subset E(b_1) \otimes B_{l}$
and $E(b_1) \otimes B_{l}$ is a finite set. Hence $B(\lam) \otimes
B_{l}$ has a head. Now, it is clear that
$u_{\lam} \otimes B_{l}$ is stable under $\te_i$'s
$(i\in I)$. Moreover, for any $u\otimes b \in B(\lam) \otimes B_{l}$,
by applying $\te_i$'s repeatedly, we get
$\te_{i_k} \ldots \te_{i_1} (u\otimes b) =u_{\lam} \otimes b'
\in u_{\lam} \otimes B_{l}$ for sufficiently large $k\ge 1$.
Hence our assertion follows from Lemma \ref{Lemma 3.5} (a).
\qed

\vskip 2mm

In the following, we will show that
the head $H(B(\lam) \otimes B_{l})$ of $B(\lam) \otimes B_{l}$
is isomorphic to the perfect crystal $B_{l-k}$. 
Moreover, we will prove that there exists an isomorphism of crystals
$$B(\lam) \otimes B_{l} \cong B_{l-k} \otimes B(\lam'),$$
where $\lam'$ is the dominant integral weight of level $k$
determined by the crystal isomorphism 
$$B(\lam) \otimes B_{k} \cong B(\lam')$$
given in \cite{KMN1}. 

\vskip 2mm

In order to give more precise statements,
let us recall the theory of  coherent families
of perfect crystals developed in \cite{KKM}.
Let $\{B_l\}_{l\ge1}$ be a family of perfect crystals $B_l$ of level $l$,
and set $B^\min_l=\{b\in B_l\,|\,\lan c,\eps(b)\ran=l\}$.
Here $\eps(b)=\sum_i\eps_i(b)\Lambda_i$, and
we will also use
$\varphi(b)=\sum_i\varphi_i(b)\Lambda_i$.
By the definition of perfect crystal, 
$\eps$ and $\varphi$ map
$B_l^\min$ bijectively to
$(P_\cl^+)_l\overset{def}{=}\{\lam\in P_\cl\mid
\lan h_i,\lam\ran\ge0,\,\lan c,\lam\ran=l\}$.
We set $J=\{(l,b)\,|\,l\ge 1,\,b\in B_l^{\min}\}$.
\Definition
\textnormal{
A crystal $B_\infty$ with an element $b_\infty$ is called 
a {\it limit} of $\{B_l\}_{l\ge 1}$
if it satisfies the following conditions:}
\eq
&&\wt(b_\infty)=0,\ \eps(b_\infty)=\varphi(b_\infty)=0,\\
&&\textnormal{for any $(l,b)\in J$, there exists an embedding of crystals}\\
&&\qquad\qquad
f_{(l,b)}:T_{\eps(b)}\otimes B_l\otimes T_{-\varphi(b)}\to B_\infty\nn\\
&&\textnormal{sending  $t_{\eps(b)}\otimes b\otimes t_{-\varphi(b)}$
to $b_\infty$,}\nn\\
&&B_\infty=\bigcup_{(l,b)\in J}\IM f_{(l,b)}.
\endeq

\end{definition}

If a limit exists, we call $\{B_l\}_{l\ge 1}$
a {\it coherent} family of perfect crystals.
It was proved in \cite{KKM} that the limit $(B_\infty,b_\infty)$
is unique up to an isomorphism.
Note that we have
\eqn
&&\lan c,\eps(b)\ran\ge0\quad\text{for any $b\in B_\infty$.}
\endeqn
We set
$B_\infty^\min=\{b\in B_\infty\,|\,\lan c,\eps(b)\ran=0\}$.
Then both $\eps$ and $\varphi$ map 
$B_\infty^\min$ 
bijectively to $P_\cl^0=\{\lam\in P_\cl\,|\,\lan c,\lam\ran=0\}$.
Moreover, there is a linear automorphism $\sigma$ of
$P_\cl^0$ such that $\sigma\varphi(b)=\eps(b)$
for any $b\in B_\infty^\min$.
We assume further the following condition:
\eq\label{cond:sigma}
&\text{$\sigma$ extends to a linear automorphism $\sigma$ of
$P_\cl$
such that}\qquad\\
&\quad\sigma\varphi(b)=\eps(b)\quad\text{for any $b\in B_l^\min$.}\nn
\endeq
We conjecture that all the coherent families satisfy this condition.
Moreover, $\sigma$ sends the simple roots to the simple roots,
and there exists an element of the Weyl group
$W$ such that its induced action on $P_\cl^0$
coincides with $\sigma|_{P_\cl^0}$.
\vskip 2mm

In the sequel, we fix a coherent family
$\{B_l\}_{l\ge 1}$ of perfect crystals 
satisfying the condition (\ref{cond:sigma}).
For positive integers $k$ and $l$ with $k<l$,
let  $\lam$ be a dominant integral weight of level $k$ and set 
$\lam'=\sigma^{-1}\lam$. Then we have:

\Lemma\label{lemma:embedding}
There exists a unique embedding of crystals
$$\psi:B_{l-k}\to T_\lam\otimes B_l\otimes T_{-\lam'}\,.$$
Moreover, we have
$\psi(B_{l-k}^\min)\subset T_\lam\otimes B_l^\min\otimes T_{-\lam'}$.
\enlemma
\proof
Let us first prove the uniqueness.
If $b\in B_{l-k}$ is sent to
$t_\lam\otimes b'\otimes t_{-\lam'}$, then
we have
$\eps(b')=\eps(b)+\lam$,
and hence we have
$\lan c,\eps(b')\ran=\lan c,\eps(b)\ran+k$.
Therefore, $\psi$ sends $B_{l-k}^\min$ to 
$T_\lam\otimes B_l^\min\otimes T_{-\lam'}$,
and $\psi|{}_{B_{l-k}^\min}$ is uniquely determined
because $\varepsilon: B_{l}^{\min} \to P_{\cl}$ is injective.
Now, the uniqueness of $\psi$ follows from the connectedness of
$B_{l-k}$.

\vskip 2mm

We shall prove the existence.
Let us take a dominant integral weight $\xi$ of level $l-k$ and
set $\mu=\lam+\xi$.
Then $\mu$ is of level $l$.
Set $\mu'=\sigma^{-1}\mu$ and $\xi'=\sigma^{-1}\xi$.
Let us take $b_l\in B_l$ such that
$\eps(b_l)=\mu$ and
$b_{l-k}\in B_{l-k}$ such that
$\eps(b_{l-k})=\xi$.
Then they are minimal vectors and we have the embeddings
\eqn
&&f_{(l,b_l)}:T_\mu\otimes B_l\otimes T_{-\mu'}\to B_\infty, \\
&&f_{(l-k,b_{l-k})}:T_\xi\otimes B_{l-k}\otimes T_{-\xi'}\to B_\infty
\endeqn
such that
$f_{(l,b_l)}(b_l)=f_{(l-k,b_{l-k})}(b_{l-k})=b_\infty$.
We shall show 
\eq\label{incl}
&&\IM (f_{(l-k,b_{l-k})})\subset\IM (f_{(l,b_l)}).
\endeq
Since $B_{l-k}$ is connected,
it is enough to show that
if $b\in B_{l-k}$ satisfies
$\te_i(b)\not=0$
and
$f_{(l-k,b_{l-k})}(t_\xi\otimes b\otimes t_{-\xi'})\in
\IM (f_{(l,b_l)})$, then
$f_{(l-k,b_{l-k})}(t_\xi\otimes \te_ib\otimes t_{-\xi'})$
also belongs to
$\IM (f_{(l,b_l)})$.
Write $f_{(l-k,b_{l-k})}(t_\xi\otimes b\otimes t_{-\xi'})
=f_{(l,b_{l})}(t_\mu\otimes b'\otimes t_{-\mu'})$
with $b'\in B_l$.
Then we have
$\eps_i(t_\xi\otimes b\otimes t_{-\xi'})
=\eps_i(t_\mu\otimes b'\otimes t_{-\mu'})$,
which implies
$\eps_i(b')=\eps_i(b)+\lan h_i,\mu-\xi\ran>0$.
Hence we have
$f_{(l-k,b_{l-k})}(t_\xi\otimes \te_ib\otimes t_{-\xi'})
=f_{(l,b_{l})}(t_\mu\otimes \te_ib'\otimes t_{-\mu'})$,
which gives (\ref{incl}). Therefore we obtain
an embedding of crystal
$T_\xi\otimes B_{l-k}\otimes T_{-\xi'}\to
T_\mu\otimes B_l\otimes T_{-\mu'}$.
This induces the desired embedding $\psi$.
\qed

\Theorem \label{theorem:iso}
Suppose  $\rank \, \Gg>2$, and let 
$\{B_l\}_{l\ge1}$ be a coherent family of perfect crystals
satisfying the condition $(\ref{cond:sigma})$.
For a pair of positive integers $k$ and $l$  with $k<l$,
let $\lam$ be a dominant integral weight of level $k$
and $\lam'=\sigma^{-1}\lam$.
Then we have an isomorphism of crystals
\begin{equation} 
B(\lam) \otimes B_{l} \cong B_{l-k} \otimes B(\lam').
\end{equation}
\entheorem
\proof
Let $\psi:B_{l-k}\to T_\lam\otimes B_l\otimes T_{-\lam'}$
be the embedding given in Lemma \ref{lemma:embedding}.
Let $B_l^{(\lam)}$ be the subset of $B_l$ such that
$\psi(B_{l-k})=T_\lam\otimes B_l^{(\lam)}\otimes T_{-\lam'}$.
In order to prove the theorem, we shall show:
\begin{equation}
\begin{align}
&H_{\lam}=u_{\lam} \otimes B_l^{(\lam)} \ \text {is closed under 
$\tilde e_{i}$'s $(i\in I)$},\label{5.3} \\
&\text {for any $b\in B_{l}$}, 
E(u_{\lam} \otimes b) \ni u_{\lam} \otimes b' \ \ \text {for some} \ 
b' \in B_l^{(\lam)}, \label{5.4} \\
&\text {there exists a bijection}  \ \ 
\Psi: u_{\lam} \otimes B_l^{(\lam)} \to B_{l-k} \label{5.5} \ \ 
  \text {that commutes} \\
& \text {with $\tilde e_{i}$'s $(i\in I)$.} \notag 
\end{align}
\end{equation}

Once we have proved them, 
Lemma \ref{Lemma 3.5} along with Lemma \ref{Lemma 3.2} would imply
$$H(B(\lam) \otimes B_{l})=u_{\lam} \otimes B_l^{(\lam)},$$
and, since $H_{\lam} \cong B_{l-k}$ is connected, 
Theorem \ref{Theorem 4.7} yields a crystal isomorphism 
$$B(\lam) \otimes B_{l} \cong H_{\lam} \otimes B(\lam')
\cong B_{l-k} \otimes B(\lam').$$

\noindent
{\sl Proof of} (\ref{5.3}) {\sl and }(\ref{5.5}):\quad
They are easily deduced from the existence of
$\psi$ and the fact
that $\te_i(b)=0$ if and only if $\eps_i(b)=0$
for $b$ in $B_{l-k}$ or in $u_\lam\otimes B_l^{(\lam)}$.

\vskip 3mm

\noindent
{\sl Proof of} (\ref{5.4}):\quad
Let us take a dominant integral weight $\xi$ of level $l-k$ and 
set $\mu=\lambda+\xi$.
Since $B_{l}$ is perfect,
there exists a unique element $b'\in B_{l}$ with $\varepsilon(b')=\mu$.
Then $b'$ belongs to $B_l^{(\lam)}$ by Lemma \ref{lemma:embedding}.
We have a crystal isomorphism 
$B(\mu) \otimes B_{l} \isomo B(\mu')$ given by
$u_{\mu} \otimes b' \mapsto u_{\mu'}$, where
$\mu'=\sigma^{-1}\mu$,
and $u_{\mu}$ (resp. $u_{\mu'}$) denotes the highest weight 
vector of $B(\mu)$ (resp. $B(\mu')$) (cf. \cite {KMN1}). 
Hence, for any $b\in B_{l}$, there exist $i_{1}, \ldots, i_{t} \in I$ 
such that 
$$\tilde e_{i_t} \ldots \tilde e_{i_1} (u_{\mu} \otimes b)
=u_{\mu} \otimes \tilde e_{i_t} \ldots \tilde e_{i_1}b
=u_{\mu} \otimes b'.$$
In particular, we have
$\varepsilon_{i_s}(\tilde e_{i_{s-1}} \ldots
\tilde e_{i_1} b)>\lan h_{i_s},\mu\ran\ge\lan h_{i_s},\lam\ran$ 
for $s=1, \ldots, t$.
This gives
$$\tilde e_{i_t} \ldots \tilde e_{i_1} (u_{\lam} \otimes b)
=u_{\lam} \otimes \tilde e_{i_t} \ldots \tilde e_{i_1}b
=u_{\lam} \otimes b' \in u_{\lam} \otimes B_l^{(\lam)},$$
which proves (\ref{5.4}).
\qed

\vskip 2mm
In the following, we will give a list of coherent families of perfect crystals
$\{B_{l}\}_{l\ge 1}$ satisfying the condition (\ref{cond:sigma})
for each quantized affine algebra 
$\Us$ of type $A_{n}^{(1)}$, $B_{n}^{(1)}$, $C_{n}^{(1)}$,
$D_{n}^{(1)}$, $A_{2n-1}^{(2)}$, $A_{2n}^{(2)}$, and 
$D_{n+1}^{(2)}$. 
For a positive integer $k<l$, and a dominant integral weight 
$\lam =a_{0} \Lam_{0} + a_{1} \Lam_{1} + \cdots
+ a_{n} \Lam_{n}$ of level $k$, Theorem \ref{theorem:iso} yields
an isomorphism of crystals
$$B(\lam) \otimes B_{l} \cong B_{l-k} \otimes B(\lam'),$$
where $\lam'=\sigma^{-1} \lam$. We will also give 
explicit descriptions of the head $u_{\lam} \otimes B_l^{(\lam)}$
of $B(\lam) \otimes B_{l}$, $\lam'=\sigma^{-1} \lam$, and the
isomorphism $\Psi: u_{\lam} \otimes B_l^{(\lam)} \isomo B_{l-k}$.
We follow the notations in \cite{KKM} and \cite{KMN2}.
\vskip 3mm

\nid
{\rm (a)} $\Gg  =A_n^{(1)} \ (n\ge 2)$: 
$$
\aligned
 B_{l}&=\{b=(x_1, \ldots, x_{n+1}) \in {\Bbb Z}_{\ge 0}^{n+1} \ 
| \ s(b)=\sum_{i=1}^{n+1} x_i =l \},\\
k&=a_0+\cdots+a_n,\\
 \lam'&=a_{n} \Lam_{0}+a_{0}\Lam_{1}+ \cdots +
a_{n-1} \Lam_{n},\\
 B_l^{(\lam)}&=\{b=(x_1, \ldots, x_{n+1}) \in B_{l} \ | 
\ x_{1} \ge a_{0}, x_{2} \ge a_{1}, \ldots, x_{n+1} \ge a_{n} \}.
\endaligned
$$
As an $A_n$-crystal, $B_l$ is isomorphic to $B(l\Lam_1)$.
The crystal structure on $B_{l}$ is described in 
\cite{KKM} and \cite{KMN2}.

\vskip 2mm

The isomorphism 
$\Psi: u_{\lam} \otimes B_l^{(\lam)} \isomo B_{l-k}$ is given by
\begin{equation}
\Psi(u_{\lam} \otimes (x_{1}, \ldots, x_{n+1}))
=(x_{1}-a_{0}, \ldots, x_{n+1}-a_n).
\end{equation}

\vskip 2mm
\nid
{\rm (b)} $\Gg =A_{2n-1}^{(2)} \ (n\ge 3)$:
$$
\aligned
 B_{l}&=\{b=(x_1, \ldots, x_{n}, \bar x_{n}, \ldots, \bar x_1) 
\in {\Bbb Z}_{\ge 0}^{2n}  | \ 
s(b)=\sum_{i=1}^{n} x_i + \sum_{i=1}^{n} \bar x_i =l \},\\
k&=a_0+a_1+2(a_2+\cdots+a_n),\\
\lam' &=a_{1} \Lam_{0}+a_{0}\Lam_{1}+ a_{2} \Lam_{2} + \cdots +
a_{n} \Lam_{n},\\
 B_l^{(\lam)}& =\{b=(x_1, \ldots, x_{n}, \bar x_{n}, \ldots, \bar x_1) 
 \in B_{l} | \\
&\qquad 
 x_i, \bar x_i \ge a_{i} \ (i=2, \ldots, n),
\ x_{1} \ge a_{0}, \ \bar x_{1} \ge a_{1}\}.
\endaligned
$$
As a $C_n$-crystal, $B_l$ is isomorphic to $B(l\Lam_1)$.
The crystal structure on $B_{l}$ is described in 
\cite{KKM} and \cite{KMN2}.

\vskip 2mm

The isomorphism 
$\Psi: u_{\lam} \otimes B_l^{(\lam)} \isomo B_{l-k}$
is given by
\begin{equation} \label{5.7}
\begin{aligned}
& \Psi(u_{\lam} \otimes (x_{1}, \ldots, x_{n},
\bar x_{n}, \ldots, \bar x_{1}))  \\
&  =(x_{1}-a_{0}, x_{2} -a_{2}, \ldots, x_{n}-a_n,
\bar x_{n} -a_{n}, \ldots, \bar x_{2}-a_{2}, \bar x_{1} -a_{1}).
\end{aligned}
\end{equation}

\vskip 2mm
\nid
{\rm (c)} $\Gg =B_n^{(1)} \ (n\ge 3)$:
$$
\aligned
 B_{l}=&\{ b=(x_1, \ldots, x_{n}, x_{0}, \bar x_{n}, \ldots, \bar x_1) 
\in {\Bbb Z}_{\ge 0}^{2n+1}
\mid\\
&\hskip 1cm  x_{0} =0 \ \text {or} \ 1, 
\ s(b)=\sum_{i=1}^{n} x_i
+x_{0} + \sum_{i=1}^{n} \bar x_i =l \},\\
k=&a_0+a_1+2(a_2+\cdots+a_{n-1})+a_n,\\
 \lam'=&a_{1} \Lam_{0}+a_{0} \Lam_{1}+ a_{2} \Lam_{2} + \cdots +
a_{n} \Lam_{n},\\
B_l^{(\lam)}=&\{b=(x_1, \ldots, x_{n}, x_0, \bar x_{n}, \ldots, \bar x_1) 
 \in B_{l} \ | \ x_{1} \ge a_{0}, \ \bar x_{1} \ge a_{1}, \\
& \hskip 1cm  x_i, \bar x_i \ge a_{i} \  (i=2, \ldots, n-1), 
\ 2x_{n}+x_0 \ge a_{n}, \ 2\bar x_n + x_{0} \ge a_{n} \}.
\endaligned
$$
As a $B_n$-crystal, $B_l$ is isomorphic to $B(l\Lam_1)$.
The crystal structure on $B_{l}$ is described in 
\cite{KKM} and \cite{KMN2}.

\vskip 2mm

The isomorphism 
$\Psi: u_{\lam} \otimes B_l^{(\lam)} \isomo B_{l-k}$ is given 
as follows. If $a_{n}$ is even, 
\begin{equation} \label{5.8}
\begin{aligned}
& \Psi(u_{\lam} \otimes (x_{1}, \ldots, x_{n}, x_{0},
\bar x_{n}, \ldots, \bar x_{1}))  \\
&  =(x_{1}-a_{0}, x_{2} -a_{2}, \ldots, x_{n}-\frac{a_n}{2}, x_{0},
\bar x_{n} -\frac{a_{n}}{2}, \ldots, \bar x_{2}-a_{2}, \bar x_{1}-a_{1}).
\end{aligned}
\end{equation}

If $a_{n}$ is odd, 
\begin{equation} \label{5.9}
\begin{aligned}
& \Psi(u_{\lam} \otimes (x_{1}, \ldots, x_{n}, x_{0},
\bar x_{n}, \ldots, \bar x_{1}))  \\
& = 
\left\{
\begin{array}{rl}
(x_{1}-a_{0}, x_{2} -a_{2}, \ldots, x_{n}-\frac{a_n +1}{2}, 1,
\bar x_{n} -\frac{a_{n}+1}{2},\ &\\
\quad \bar x_{n-1} - a_{n-1},  \ldots, \bar x_{2}-a_{2}, \bar x_{1} -a_{1})
& \text {if} \ x_{0}=0, \\
 (x_{1}-a_{0}, x_{2} -a_{2}, \ldots, x_{n}-\frac{a_n -1}{2}, 0,
 \bar x_{n} -\frac{a_{n}-1}{2},\ & \\
\quad \bar x_{n-1} - a_{n-1}, \ldots, \bar x_{2}-a_{2}, \bar x_{1} -a_{1})
& \text {if} \ x_{0}=1.
\end{array}\right.
\end{aligned}
\end{equation}

\vskip 2mm
\nid
{\rm (d)} $\Gg=A_{2n}^{(2)} \ (n\ge 2)$:
$$
\aligned
 B_{l}=&\{b=(x_1, \ldots, x_{n}, \bar x_{n}, \ldots, \bar x_1) 
\in {\Bbb Z}_{\ge 0}^{2n} \ | \ 
s(b)= \sum_{i=1}^{n} x_i +\sum_{i=1}^{n} \bar x_i \le l \},\\
k=&a_0+2(a_1+\cdots+a_n),\\
 \lam'=&\lam=a_{0} \Lam_{0}+a_{1}\Lam_{1} + \cdots +
a_{n} \Lam_{n},\\
 B_l^{(\lam)}=&\{b=(x_1, \ldots, x_{n}, \bar x_{n}, \ldots, \bar x_1) 
 \in B_{l} \mid
 x_i, \bar x_i \ge a_{i} \  (i=1, \ldots, n), 
\  s(b) \le l-a_{0} \}.
\endaligned
$$
As a $C_n$-crystal, $B_l$ is isomorphic to 
$B(0) \oplus B(\Lam_{1}) \oplus \cdots \oplus B(l\Lam_1)$.
The crystal structure on $B_{l}$ is described in 
\cite{KKM} and \cite{KMN2}.

\vskip 2mm

The isomorphism 
$\Psi: u_{\lam} \otimes B_l^{(\lam)} \isomo B_{l-k}$
is given by
\begin{equation} 
\begin{aligned}
& \Psi(u_{\lam} \otimes (x_{1}, \ldots, x_{n},
\bar x_{n}, \ldots, \bar x_{1}))  \\
&  =(x_{1}-a_{1}, x_{2} -a_{2}, \ldots, x_{n}-a_n,
\bar x_{n} -a_{n}, \ldots, \bar x_{2}-a_{2}, \bar x_{1} -a_{1}).
\end{aligned}
\end{equation}

\vskip 2mm
\nid
{\rm (e)} $\Gg=D_{n+1}^{(2)} \ (n\ge 2)$:
$$
\aligned
B_{l}=&\{ b=(x_1, \ldots, x_{n}, x_0, \bar x_{n}, \ldots, \bar x_1) 
\in {\Bbb Z}_{\ge 0}^{2n+1}\mid\\
&
\hskip 1cm 
x_{0}=0 \ \text {or} \ 1,
\ s(b)= \sum_{i=1}^{n} x_i +x_{0}+ \sum_{i=1}^{n} \bar x_i \le l \},\\
k=&a_0+2(a_1+\cdots+a_{n-1})+a_n,\\
\lam'=&\lam=a_{0} \Lam_{0}+a_{1}\Lam_{1} + \cdots +
a_{n} \Lam_{n},\\
B_l^{(\lam)}=&\{ b=(x_1, \ldots, x_{n}, x_{0}, \bar x_{n}, \ldots, \bar x_1) 
\in B_{l} \mid\ x_i, \bar x_i \ge a_{i} \ (i=1, \ldots, n-1),\\
&
\hskip 1cm \ 2x_n+x_0 \ge a_n, \ 2\bar x_n + x_0 \ge a_n, 
\ s(b) \le l-a_{0} \}.
\endaligned
$$
As a $B_n$-crystal, $B_l$ is isomorphic to 
$B(0) \oplus B(\Lam_{1}) \oplus \cdots \oplus B(l\Lam_1)$.
The crystal structure on $B_{l}$ is described in 
\cite{KKM} and \cite{KMN2}.

\vskip 2mm

The isomorphism 
$\Psi: u_{\lam} \otimes B_l^{(\lam)} \isomo B_{l-k}$ is given 
as follows. If $a_{n}$ is even, 
\begin{equation} 
\begin{aligned}
& \Psi(u_{\lam} \otimes (x_{1}, \ldots, x_{n}, x_{0},
\bar x_{n}, \ldots, \bar x_{1}))  \\
&\quad =(x_{1}-a_{1}, x_{2} -a_{2}, \ldots, x_{n}-\frac{a_n}{2}, x_{0},
\bar x_{n} -\frac{a_{n}}{2}, \ldots, \bar x_{2}-a_{2}, \bar x_{1} -a_{1}).
\end{aligned}
\end{equation}
If $a_{n}$ is odd, 
\begin{equation} 
\begin{aligned}
& \Psi(u_{\lam} \otimes (x_{1}, \ldots, x_{n}, x_{0},
\bar x_{n}, \ldots, \bar x_{1}))  \\
& = 
\left\{
\begin{array}{rl}
(x_{1}-a_{1}, x_{2} -a_{2}, \ldots, x_{n}-\frac{a_n +1}{2}, 1,
\bar x_{n} -\frac{a_{n}+1}{2},\ &\\
\quad \bar x_{n-1} - a_{n-1},  \ldots, \bar x_{2}-a_{2}, \bar x_{1} -a_{1})
& \text {if} \ x_{0}=0, \\
 (x_{1}-a_{1}, x_{2} -a_{2}, \ldots, x_{n}-\frac{a_n -1}{2}, 0,
 \bar x_{n} -\frac{a_{n}-1}{2},\ & \\
\quad \bar x_{n-1} - a_{n-1}, \ldots, \bar x_{2}-a_{2}, \bar x_{1} -a_{1})
& \text {if} \ x_{0}=1.
\end{array}\right.
\end{aligned}
\end{equation}

\vskip 2mm
\nid
{\rm (f)} $\Gg=C_n^{(1)} \ (n\ge 2)$:
$$
\aligned
B_{l}=&\{b=(x_1, \ldots, x_{n}, \bar x_{n}, \ldots, \bar x_1) 
\in {\Bbb Z}_{\ge 0}^{2n} \mid
s(b)= \sum_{i=1}^{n} x_i +\sum_{i=1}^{n} \bar x_i \le 2l,
\ s(b) \in 2 {\Bbb Z} \},\\
k=&a_0+\cdots+a_n,\\
\lam'=&\lam=a_{0} \Lam_{0}+a_{1}\Lam_{1} + \cdots +
a_{n} \Lam_{n},\\
B_l^{(\lam)}=&\{b=(x_1, \ldots, x_{n}, \bar x_{n}, \ldots, \bar x_1) 
\in B_{l} \mid
x_i,\bar x_i\ge a_{i}\ (i=1, \ldots, n),
\ s(b) \le 2(l-a_{0}) \}.
\endaligned
$$
As a $C_n$-crystal, $B_l$ is isomorphic to 
$B(0) \oplus B(2\Lam_{1}) \oplus \cdots \oplus B(2l\Lam_1)$.
The crystal structure on $B_{l}$ is described in 
\cite{KKM}.

\vskip 2mm

The isomorphism 
$\Psi: u_{\lam} \otimes B_l^{(\lam)} \isomo B_{l-k}$
is given by
\begin{equation} 
\begin{aligned}
& \Psi(u_{\lam} \otimes (x_{1}, \ldots, x_{n},
\bar x_{n}, \ldots, \bar x_{1}))  \\
&\quad  =(x_{1}-a_{1}, x_{2} -a_{2}, \ldots, x_{n}-a_n,
\bar x_{n} -a_{n}, \ldots, \bar x_{2}-a_{2}, \bar x_{1} -a_{1}).
\end{aligned}
\end{equation}

\vskip 2mm
\nid
{\rm (g)} $\Gg=D_n^{(1)} \ (n\ge 4)$:
$$
\aligned
B_{l}=&\{b=(x_1, \ldots, x_{n}, \bar x_{n}, \ldots, \bar x_1) 
\in {\Bbb Z}_{\ge 0}^{2n} \mid
x_n =0 \ \text {or} \ \bar x_n=0,
\  s(b)=\sum_{i=1}^{n} x_i
+ \sum_{i=1}^{n} \bar x_i =l \},\\
k=&a_0+a_1+2(a_2+\cdots+a_{n-2})+a_{n-1}+a_n,\\
\lam'=& a_{1} \Lam_{0}+a_{0}\Lam_{1}+ a_{2} \Lam_{2} + \cdots +
a_{n-2} \Lam_{n-2} +a_n \Lam_{n-1}+a_{n-1} \Lam_{n},\\
B_l^{(\lam)} =& 
\left\{
\begin{array}{ll}
 \{b=(x_1, \ldots, x_{n}, \bar x_{n}, \ldots, \bar x_1) 
 \in B_{l} \ | \ x_{1} \ge a_{0}, \ \bar x_{1} \ge a_{1},&\\
\hskip 1cm x_i, \bar x_i \ge a_{i} \ (i=2, \ldots, n-2),
\ x_{n-1}, \bar x_{n-1} \ge a_{n},&\\
\hskip 1cm x_{n-1}+x_{n} \ge a_{n-1}, \bar x_{n-1}+x_{n} \ge a_{n-1} \} 
&  \text {if} \ a_{n-1} \ge a_{n}, \\
 \{b=(x_1, \ldots, x_{n}, \bar x_{n}, \ldots, \bar x_1) 
 \in B_{l} \ | \  x_{1} \ge a_{0}, \ \bar x_{1} \ge a_{1}, &\\
 \hskip 1cm x_i, \bar x_i \ge a_{i} \ (i=2, \ldots, n-2), 
\ x_{n-1}, \bar x_{n-1} \ge a_{n-1}, &\\
 \hskip 1cm 
x_{n-1}+\bar x_{n} \ge a_{n}, \ \bar x_{n-1}+ \bar x_{n} \ge a_{n} \}
& \text {if} \ a_{n-1} \le a_{n}.
\end{array}\right.
\endaligned
$$
As a $D_n$-crystal, $B_l$ is isomorphic to $B(l\Lam_1)$.
The crystal structure on $B_{l}$ is described in 
\cite{KKM} and \cite{KMN2}.

\vskip 2mm

The isomorphism $\Psi: u_{\lam} \otimes B_l^{(\lam)}\isomo B_{l-k}$ is given 
as follows. If $a_{n-1} \ge a_{n}$,
\begin{equation}
\begin{aligned}
& \Psi(u_{\lam} \otimes (x_{1}, \ldots, x_{n},
\bar x_{n}, \ldots, \bar x_{1}))  \\
& \hskip 3mm  =(x_{1}-a_{0}, x_{2} -a_{2}, \ldots, x_{n-2}-a_{n-2}, \\
& \hskip 1cm  x_{n-1}-a_{n} -(a_{n-1}-a_{n}-x_{n})_{+},
\ (x_{n}-a_{n-1}+a_{n})_{+},  \\
& \hskip 1cm \bar x_n + (a_{n-1}-a_{n}-x_{n})_{+}, 
\ \bar x_{n-1}-a_{n}-(a_{n-1}-a_{n}-x_{n})_{+}, \\
& \hskip 1cm 
 \bar x_{n-2} -a_{n-2}, \ldots, \bar x_{2}-a_{2}, \bar x_{1} -a_{1}),
\end{aligned}
\end{equation}
and if $a_{n-1} \le a_{n}$,
\begin{equation}
\begin{aligned}
& \Psi(u_{\lam} \otimes (x_{1}, \ldots, x_{n},
\bar x_{n}, \ldots, \bar x_{1}))  \\
& \hskip 3mm  =(x_{1}-a_{0}, x_{2} -a_{2}, \ldots, x_{n-2}-a_{n-2}, \\
& \hskip 1cm x_{n-1}-a_{n-1}-(a_{n}-a_{n-1}-\bar x_{n})_{+},
\ x_n+(a_{n}-a_{n-1}- \bar x_{n})_{+}, \\
& \hskip 1cm  (\bar x_{n}-a_{n}+a_{n-1})_{+}, 
\ \bar x_{n-1}-a_{n-1} -(a_{n}-a_{n-1}- \bar x_{n})_{+}, \\
& \hskip 1cm
 \bar x_{n-2} -a_{n-2}, \ldots, \bar x_{2}-a_{2}, \bar x_{1} -a_{1}).
\end{aligned}
\end{equation}

%


\bibliographystyle{unsrt}

\begin{thebibliography}{10}

\bibitem{Akas}
T. Akasaka, M. Kashiwara, 
\newblock{\it Finite-dimensional representations 
of quantum affine algebras},
\newblock{to appear in Publ. RIMS, q-alg\thinspace9703028}.

\bibitem{Kac}
V. Kac,
\newblock{\it Infinite Dimensional Lie Algebras}, 3rd ed., 
\newblock{Cambridge University Press}, 
\newblock 1990.

\bibitem{Kas1}
M. Kashiwara,
\newblock{\it On crystal bases of the $q$-analogue of universal
enveloping algebras},
\newblock{Duke Math. J.} {\bf 63} 
\newblock (1991), 465--516.

\bibitem{Kas2}
M. Kashiwara, 
\newblock{\it Crystal bases of modified quantized
enveloping algebras},
\newblock{Duke Math. J.} {\bf 73} 
\newblock (1994), 383--413.

\bibitem{KKM}
S.-J. Kang, M. Kashiwara, K. C. Misra, 
\newblock{\it Crystal bases of Verma modules 
for quantum affine Lie algebras},
\newblock{Compositio Math.} {\bf 92} 
\newblock (1994), 299--325.

\bibitem{KMN1}
S.-J. Kang, M. Kashiwara, K. C. Misra, T. Miwa, 
T. Nakashima, A. Nakayashiki, 
\newblock{\it Affine crystals and vertex models},
\newblock{Int. J. Mod. Phys.} {\bf A, Suppl. 1A} 
\newblock (1992), 449--484.

\bibitem{KMN2}
\bysame,
\newblock{\it Perfect crystals for quantum affine Lie algebras},
\newblock{Duke Math. J.} {\bf 68} 
\newblock (1992), 499--607.

\bibitem{Lusztig}
G. Lusztig,
\newblock{\it Introduction to Quantum Groups}, 
Progress of Mathematics {\bf 10},
\newblock{Birkh\"auser}, 
\newblock 1993.

\bibitem{Nak}
A. Nakayashiki,
\newblock{\it Fusion of the $q$-vertex operators and its applications to
solvable vertex models}, 
\newblock{Commun. Math. Phys.} {\bf 177}
\newblock (1996), 27-62.

\bibitem{Saito}
Y. Saito,
\newblock{\it PBW basis of quantized universal enveloping algebras}, 
\newblock{Publ. RIMS} {\bf 30}
\newblock (1994), 209--232.

\end{thebibliography}

\end{document}